\documentclass[a4paper,11pt]{article}
\usepackage{latexsym}
\usepackage{amssymb}
\usepackage{amsfonts}
\usepackage{amsmath}
\usepackage{indentfirst}
\usepackage{graphicx}
\usepackage{epsfig}
\usepackage[enableskew]{youngtab}

\usepackage{fullpage}

\usepackage{color}
\usepackage{mathrsfs}
\usepackage{manfnt}

\usepackage{tikz,ifthen}

\newcounter{indice}

\newcommand{\shape}[1]{
\setcounter{indice}{0};
\foreach \i in {#1} {
\addtocounter{indice}{1};
\foreach \x in {1,...,\i} {
\draw (\x-1,-\theindice+1) rectangle (\x,-\theindice);
}
}
}

\newcommand{\griddedshape}[1]{
\setcounter{indice}{0};
\foreach \i in {#1} {
\addtocounter{indice}{1};
\foreach \x in {1,...,\i} {
\draw[gray,dotted] (\x-1,-\theindice+1) rectangle (\x,-\theindice);
}
}
}

\newcommand{\skewshape}[1]{
\setcounter{indice}{0};
\foreach \i/\j in {#1} {
\addtocounter{indice}{1};
\foreach \x in {\j,...,\i} {
\draw (\x-1,-\theindice+1) rectangle (\x,-\theindice);
}
}
}

\newtheorem{prop}{Proposition}[section]
\newtheorem{teor}{Theorem}[section]

\newcommand{\cvd}{\hfill $\blacksquare$\bigskip}
\date{}

\author{Luca Ferrari\thanks{Dipartimento di Sistemi e Informatica, viale Morgagni 65, 50134 Firenze, Italy
{\tt ferrari@dsi.unifi.it}} \and Emanuele
Munarini\thanks{Politecnico di Milano, Dipartimento di Matematica,
Piazza Leonardo da Vinci 32, 20133 Milano, Italy {\tt
emanuele.munarini@polimi.it}}}

\title{Enumeration of saturated chains in Dyck lattices}

\begin{document}

\maketitle

\begin{abstract}
We determine a general formula to compute the number of saturated
chains in Dyck lattices, and we apply it to find the number of
saturated chains of length 2 and 3. We also compute what we call
the Hasse index (of order 2 and 3) of Dyck lattices, which is the
ratio between the total number of saturated chains (of length 2
and 3) and the cardinality of the underlying poset.
\end{abstract}

\section{Introduction}

Given a poset, a very natural problem is to count how many
saturated chains it has. A \emph{saturated chain} in a poset is a
chain such that, if $x<y$ are consecutive elements in the chain,
then $y$ covers $x$. In the present paper we wish to address this
problem in the case of Dyck lattices. The \emph{Dyck lattice of
order $n$}, to be denoted $\mathcal{D}_n$, is the lattice of Dyck
paths of semilength $n$ whose associated partial order relation is
given by \emph{containment}: given $\gamma,\gamma '\in
\mathcal{D}_n$, it is $\gamma \leq \gamma '$ when, in the usual
two-dimensional drawing of Dyck paths, $\gamma$ lies weakly below
$\gamma '$. Some papers studying properties of Dyck lattices are
\cite{FM1,FP}. Counting saturated chains of length 1 is clearly
equivalent to enumerating edges in the associated Hasse diagram,
which has been considered in \cite{FM2} not only for Dyck lattices
but also for other lattices of paths. Here we start by providing a
general formula for counting saturated chains of length $h$, for
any fixed $h$, in a given Dyck lattice. Next we deal with the
cases $h=2,3$, giving for them detailed enumerative results. We
also define the notion of \emph{Hasse index of order $h$} (thus
generalizing the concept of Hasse index proposed in \cite{FM2})
and compute such an index in the two mentioned special cases.

\section{Preliminaries}

In this section we collect some notations and results which will
be used in the sequel.

\bigskip

A \emph{Young tableau} is a filling of a Ferrers shape $\lambda$
using distinct positive integers from $1$ to $n=| \lambda |$, with
the properties that the values are (strictly) decreasing along
each row and each column of the Ferrers shape.
Here $|\lambda |$ denotes the number of cells of the Ferrers shape $\lambda$.
This constitutes a slight departure from the classical definition,
which requires the word ``increasing" instead of the word ``decreasing".
However, it is clear that
all the properties and results on (classical) Young tableaux
can be translated into our setting by simply replacing
the total order ``$\leq$" with the total order ``$\geq$" on $\mathbb{N}$.
A \emph{skew Young tableau} is defined exactly as a Young tableau,
with the only difference that the underlying shape
consists of a Ferrers shape $\lambda$
with a (possibly empty) Ferrers shape $\mu$ removed (starting from the top-left corner),
in such a way that the resulting shape is \emph{strongly connected}:
this means that every pair of consecutive rows has at
least one common column and every pair of consecutive columns has
at least one common row (such a shape will be also called a
\emph{skew Ferrers shape}).

A \emph{Dyck path} is a path starting from the origin of a fixed
Cartesian coordinate system, ending on the $x$-axis, never going
below the $x$-axis, and using only the two steps $u=(1,1)$ and
$d=(1,-1)$. A \emph{valley} (\emph{peak}) in a Dyck path is a pair
of consecutive steps $du$ ($ud$). The semilength of a Dyck path is
just half the number of its steps. The set of all Dyck paths of
semilength $n$ will be denoted $D_n$.

The set $D_n$ endowed with the partial order described in the
Introduction will be called the \emph{Dyck lattice of order $n$}
and denoted $\mathcal{D}_n$. The generating series of saturated
chains of length $h$ in the family of Dyck lattices will be
written $SC_h (x)$, whereas the number of saturated chains of
length $h$ in $\mathcal{D}_n$ (i.e. the coefficient of $x^n$ in
$SC_h (x)$) will be written $sc_h (\mathcal{D}_n)$.

At the end of this section, we propose a generalization of the
notion of Hasse index given in \cite{FM2}. Recall that the
\emph{Hasse index} $i(\mathcal{P})$ of a poset $\mathcal{P}$ is
given by $i(\mathcal{P})=\frac{\ell
(\mathcal{P})}{|\mathcal{P}|}$, where $\ell (\mathcal{P})$ is the
number of covering pairs in $\mathcal{P}$. Given a positive
integer $h$, we now define the \emph{Hasse index of order $h$} of
$\mathcal{P}$ as $i_h (\mathcal{P})=\frac{sc_h
(\mathcal{P})}{|\mathcal{P}|}$, where $sc_h (\mathcal{P})$ denotes
the number of saturated chains of length $h$ of the poset
$\mathcal{P}$. Of course $i_1 (\mathcal{P})=i(\mathcal{P})$. For
instance, for the Boolean algebra $\mathcal{B}_n$ having $2^n$
elements, $sc_h (\mathcal{B}_n )$ can be computed by taking an
arbitrary subset having $k$ elements (for $0\leq k\leq n$) and
then adding any $h$ of the remaining elements in a specified order.
Equivalently, we can choose a subset having $h$ elements, a linear order on it,
and a subset of its complement. Therefore we get
\begin{displaymath}
sc_h (\mathcal{B}_n )=\sum_{k=0}^{n}{n\choose k}(n-k)_h =(n)_h
2^{n-h},
\end{displaymath}
where $(a)_b =a\cdot (a-1)\cdot \ldots \cdot (a-b+1)$ denotes a \emph{falling factorial}.

Thus, the Hasse index of order $h$ of $\mathcal{B}_n$ is given by
\begin{displaymath}
i_h (\mathcal{B}_n )=\frac{(n)_h\cdot
2^{n-h}}{2^n}=\frac{(n)_h}{2^h}.
\end{displaymath}

We will say that the Hasse index of order $h$ of a sequence
$\mathscr{P}=\{ \mathcal{P}_0 ,\mathcal{P}_1 ,\mathcal{P}_2
,\ldots ,\mathcal{P}_n ,\ldots \}$ of posets is \emph{Boolean}
when $i_h (\mathcal{P}_n )=\frac{(n)_h}{2^h}$ and
\emph{asymptotically Boolean} when $i_h (\mathcal{P}_n )\sim
\frac{(n)_k}{2^h}$ (or, which is the same, $i_h (\mathcal{P}_n
)\sim \frac{n^h}{2^h}$).
%
%

In the computation of the Hasse index we will also use the well
known \emph{Darboux theorem} (see for instance \cite{BLL}), which
asserts that, given a complex number $\xi \ne 0$ and a complex
function $f(x)$ analytic at the origin, if $f(x)=(1-x/\xi
)^{-\alpha }\psi(x)$, where $\psi(x)$ is a series with radius of
convergence $R>|\xi |$ and $\alpha \not \in \{ 0,-1,-2,\ldots \}$,
then
\begin{displaymath}
[x^n]f(x)\sim \frac{\psi(\xi )}{\xi^n }\frac{n^{\alpha -1}}{\Gamma
(\alpha )}.
\end{displaymath}

\section{The general enumeration formula}

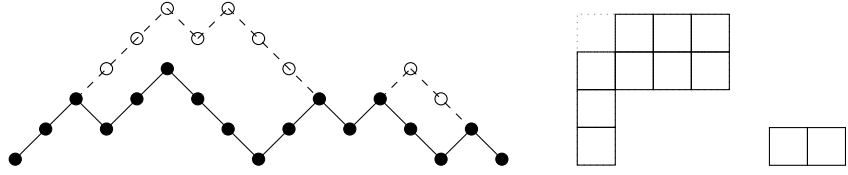
\begin{figure}[ht]
\begin{center}
\begin{tikzpicture}
\begin{scope}[scale=0.4]
\draw (0,0) -- (2,2);
\draw (2,2) -- (3,1);
\draw (3,1) -- (5,3);
\draw (5,3) -- (8,0);
\draw (8,0) -- (10,2);
\draw (10,2) -- (11,1);
\draw (11,1) -- (12,2);
\draw (12,2) -- (14,0);
\draw (14,0) -- (15,1);
\draw (15,1) -- (16,0);
\draw (0,0) [fill] circle (.2);
\draw (1,1) [fill] circle (.2);
\draw (2,2) [fill] circle (.2);
\draw (3,1) [fill] circle (.2);
\draw (4,2) [fill] circle (.2);
\draw (5,3) [fill] circle (.2);
\draw (6,2) [fill] circle (.2);
\draw (7,1) [fill] circle (.2);
\draw (8,0) [fill] circle (.2);
\draw (9,1) [fill] circle (.2);
\draw (10,2) [fill] circle (.2);
\draw (11,1) [fill] circle (.2);
\draw (12,2) [fill] circle (.2);
\draw (13,1) [fill] circle (.2);
\draw (14,0) [fill] circle (.2);
\draw (15,1) [fill] circle (.2);
\draw (16,0) [fill] circle (.2);
\draw[dashed] (2,2) -- (5,5);
\draw[dashed] (5,5) -- (6,4);
\draw[dashed] (6,4) -- (7,5);
\draw[dashed] (7,5) -- (10,2);
\draw[dashed] (12,2) -- (13,3);
\draw[dashed] (13,3) -- (15,1);
\draw (3,3) circle (.2);
\draw (4,4) circle (.2);
\draw (5,5) circle (.2);
\draw (6,4) circle (.2);
\draw (7,5) circle (.2);
\draw (8,4) circle (.2);
\draw (9,3) circle (.2);
\draw (13,3) circle (.2);
\draw (14,2) circle (.2);
\end{scope}
\end{tikzpicture} \qquad \begin{tikzpicture}
\begin{scope}[scale=0.5]
\griddedshape{4,4,1,1};
\skewshape{4/2,4/1,1/1,1/1};
\end{scope}
\end{tikzpicture} \quad \begin{tikzpicture}
\begin{scope}[scale=0.5]
\shape{2};
\end{scope}
\end{tikzpicture}
\end{center}
\caption{A pair of Dyck paths $\gamma$ (thick) and $\gamma'$ (dashed), with $\gamma <\gamma'$, and the corresponding set of skew Ferrers shapes.\label{pair}}
\end{figure}

Let $\gamma^{(0)}<\gamma^{(1)}<\cdots <\gamma^{(h)}$
be a saturated chain (of length $h$) in $\mathcal{D}_n$. It is
easy to see that two consecutive paths of the chain only differ by
a pair of consecutive steps, namely a valley (a peak) in
the smallest (largest) one. More generally, the minimum
$\gamma^{(0)}$ and the maximum $\gamma^{(h)}$ differ by a set of
steps in such a way that the sum of the areas of the regions
delimited by these steps is equal to $h$. To be more precise, this
means that the two paths can be factorized as $\gamma^{(0)}
=\alpha_1 \gamma^{(0)}_1 \alpha_2 \gamma^{(0)}_2 \cdots \alpha_k
\gamma^{(0)}_k \alpha_{k+1}$ and $\gamma^{(h)} =\alpha_1
\gamma^{(h)}_1 \alpha_2 \gamma^{(h)}_2 \cdots \alpha_k
\gamma^{(h)}_k \alpha_{k+1}$, where, for every $i$, the two
factors $\gamma^{(0)}_i$ and $\gamma^{(h)}_i$ have the same
length, and the sum of the areas of the regions determined by the
pairs of factors $(\gamma^{(0)}_i ,\gamma^{(h)}_i)$ is equal to
$h$ (see Figure \ref{pair}).

Each of the regions determined by the pairs $(\gamma^{(0)}_i
,\gamma^{(h)}_i)$ can be regarded as a \emph{skew Ferrers shape}.
To fix notations, we will suppose that such a shape is that
obtained by rotating the sheet of paper by $45^{\circ}$
anticlockwise. Referring again to Figure \ref{pair}, the pair of
Dyck paths on the left determines the pair of skew Ferrers shapes
on the right.

Now suppose to select a saturated chain from $\gamma^{(0)}$ to
$\gamma^{(h)}$. This corresponds to choosing, one at a time, the
$h$ cells belonging to the skew Ferrers shapes described above.
More formally, this defines a linear order on the set of all cells
of the skew Ferrers shapes determined by the two paths such that,
on each row and on each column, cells are in decreasing order.
This means that a saturated chain essentially generates a set of
\emph{skew Young tableaux}.

\bigskip

Let now $\gamma \in \mathcal{D}_n$. We want to determine the
number of saturated chains of length $h$ starting from $\gamma$ in
$\mathcal{D}_n$. According to the above considerations, to
describe any such chain we start by giving a partition $\lambda
=(\lambda_1 , \ldots ,\lambda_k )$ of $h$. Next we have to choose
a set $\gamma_1 ,\ldots ,\gamma_k$ of factors of $\gamma$ such
that, for any $i\leq k$, we can build a skew Ferrers shape
$\varphi_i$ on $\gamma_i$ having area $\lambda_i$. Finally, to
determine the saturated chain, we just have to linearly order the
cells of the Ferrers shapes thus obtained, or, equivalently, to
endow each of the shapes with a skew Young tableaux structure.

Now we will try to describe more formally the above argument.
Denote by $SkFS$ the set of all skew Ferrers shapes. Given
$\varphi \in SkFS$, we write $A(\varphi )$ for the \emph{area} of
$\varphi$, i.e. the number of cells of $\varphi$. We also define
$SkFS(n)=\{ \varphi \in SkFS \; |\; A(\varphi )=n\}$. Given a set
of words $\gamma_1 ,\ldots ,\gamma_n$ on the alphabet $\{ u,d\}$,
we say that they are a set of \emph{pairwise disjoint occurrences
(p.d.o.)} in $\gamma$ when they appear as factors of $\gamma$
having no pairwise intersection. A skew Ferrers shape $\varphi$ is
delimited by two paths, both starting at its bottom left corner
and ending at its top right corner. Each of such paths can be seen
as a word on $\{ u,d\}$, by simply encoding a horizontal step with
the letter $d$ and a vertical step with the letter $u$. The word
having $d$ as its first letter is called the \emph{lower border}
of $\varphi$ and is denoted $b(\varphi )$. Finally, for any given
$\varphi \in SkFS$, let $t(\varphi )$ be the number of skew Young
tableaux of shape $\varphi$.

For any path $\gamma \in \mathcal{D}_n$, let $\lambda =(\lambda_1 ,\ldots
,\lambda_k )$ be a partition of the positive integer $h$ (this
will also be written as $\lambda \vdash h$). Next we have to
choose a set $\gamma_1 ,\ldots ,\gamma_k$ of pairwise disjoint
occurrences in $\gamma$ such that, for any $i\leq k$, there exists
a skew Ferrers shape $\varphi_i \in SkFS(\lambda_i )$ for which
$b(\varphi_i )=\gamma_i$. Now, to get a saturated chain, we have
to select a $k$-tuple $(\varphi_1 ,\ldots ,\varphi_k )\in SkFS^k$
such that $b(\varphi_i )=\gamma_i$ and $A(\varphi_i )=\lambda_i$,
and for each component $\varphi_i$ we have to choose one among the
$t(\varphi_i )$ possible skew Young tableaux. Finally, since the
set of integers actually used to fill in the cells of each
$\varphi_i$ can be any possible set of $|\lambda_i |$ integers
less than or equal to $h$, we have proved the following result.

\begin{teor} The number $sc_h (\mathcal{D}_n)$ of saturated chains of
length $h$ of the lattice $\mathcal{D}_n$ is given by the
following formula:
\begin{equation}\label{general}
\sum_{\gamma \in \mathcal{D}_n}\sum_{\lambda \vdash
h}\sum_{\gamma_1 ,\ldots ,\gamma_k \textnormal{p.d.o.}\atop
(\forall i )(\exists \varphi_i \in SkFS(\lambda_i ))b(\varphi_i )
=\gamma_i}\sum_{(\varphi_1 ,\ldots ,\varphi_k )\in SkFS^k \atop
(\forall i)(b(\varphi_i )=\gamma_i ,A(\varphi_i )=\lambda_i
)}{h\choose A(\varphi_1 ),\ldots ,A(\varphi_k )}t(\varphi_1 )\cdots
t(\varphi_k ) .
\end{equation}
\end{teor}

In the rest of the paper, our main aim is to apply the above
formula to the special cases $h=2$ and $h=3$ (the case $h=1$
having already been examined in \cite{FM2}), thus finding some new
results on the poset structure of Dyck lattices.

We end the present section by recalling that this problem could
also be tackled from a slightly different point of view. Indeed,
given two Dyck paths of the same length $\gamma_1$ and $\gamma_2$
such that $\gamma_1 \leq \gamma_2$, the set of all saturated
chains between $\gamma_1$ and $\gamma_2$ can be represented by
means of a suitable \emph{Polya festoon}, more precisely a Polya
festoon whose components cannot be $-$polygons (see \cite{F}). It
seems that this approach could be more elegant, but should lead to
more difficult computations.

We also remark that, in the paper \cite{CPQS}, pairs of
noncrossing free Dyck paths (also called Grand-Dyck paths in
different sources) are considered, also in connection with several
different combinatorial structures, such as noncrossing partitions
and vacillating tableaux. It could be of some interest to extend
our results to the case of free Dyck paths and successively
interpret them on the above mentioned combinatorial objects via
the bijections described in \cite{CPQS}.


\section{Saturated chains of length 2}

In order to apply formula (\ref{general}) to the case of saturated
chains of length 2 we simply have to set $h=2$. Doing this way,
one immediately observes that there are only two partitions of 2,
namely (1,1) and (2), and that there exists one pair of ``admissible"
skew Ferrers shapes of area 1, i.e. $({\tiny \yng(1)},{\tiny
\yng(1)})$, and two different skew Ferrers shapes of area 2, i.e.
${\tiny \yng(2)}$ and ${\tiny \yng(1,1)}$ . Since each of these
shapes can be endowed with only one Young tableau structure, we
arrive at the following result.

\begin{prop} The generating series for the number of saturated
chains of length 2 of Dyck lattices is given by
\begin{equation}\label{2sat}
SC_2 (x)=\sum_{n\geq 0}\left( \sum_{\gamma \in
\mathcal{D}_n}\left( 2\cdot \# (du,du)_{\gamma}+\#
(ddu)_{\gamma}+\# (duu)_{\gamma} \right) \right) x^n ,
\end{equation}
where with $\# (\gamma_1 ,\ldots ,\gamma_k )_{\gamma}$ we denote
the number of pairwise disjoint occurrences of the $\gamma_i$'s in
$\gamma$.
\end{prop}

All we have to do now is to evaluate the three unknown
quantities appearing in (\ref{2sat}). The following proposition
translates formula (\ref{2sat}) into an expression more suitable
for computing.

\begin{prop} Denote with $F(q,x)$ and $V(q,x)$ the generating
series of all Dyck paths where $x$ keeps track of the semilength
and $q$ keeps track of the factor $duu$ and of the factor $du$
(i.e. valleys), respectively. Then
\begin{equation}\label{2genser}
SC_2 (x)=2\cdot \left[ \frac{\partial F}{\partial
q}\right]_{q=1}+\left[ \frac{\partial^2 V}{\partial
q^2}\right]_{q=1}.
\end{equation}
\end{prop}

\emph{Proof.}\quad The expression $\left[ \frac{\partial
V}{\partial q}\right]_{q=1}$ gives the generating series of Dyck
paths with respect to the number of valleys. Analogously, the
expression $\left[ \frac{\partial^2 V}{\partial q^2}\right]_{q=1}$
gives the generating series of Dyck paths with respect to the
number of (non ordered) pairs of valleys. Moreover, the expression
$\left[ \frac{\partial F}{\partial q}\right]_{q=1}$ gives the
generating series of Dyck paths with respect to the number of
factors $duu$. Since the factors $ddu$ and $duu$ are obviously
equidistributed in the set of Dyck paths, formula (\ref{2genser})
immediately follows.\cvd

We are now in a position to find a neat expression for the
generating series $SC_2 (x)$.

\begin{teor}\label{2teor} The generating series for the number of saturated
chains of length 2 of Dyck lattices is given by
\begin{equation}\label{2chains}
SC_2 (x)=\sum_{n\geq 0}sc_2 (\mathcal{D}_n )x^n =\frac{1-6x+6x^2-(1-4x)\sqrt{1-4x}}{-(1-4x)\sqrt{1-4x}},
\end{equation}
where
\begin{equation}\label{2coeff}
sc_2 (\mathcal{D}_n)={2n\choose n}\frac{(n-1)(n-2)}{2(2n-1)}\qquad (n\geq 1).
\end{equation}
\end{teor}

\emph{Proof.}\quad Let $G(q,x)$, $H(q,x)$ be the generating series of
Dyck paths starting with a peak and Dyck paths starting with two
consecutive up steps, respectively, where $x$ keeps track of the
semilength and $q$ keeps track of the factor $duu$. Since any
non-empty Dyck path $\gamma$ decomposes uniquely as
$\gamma=U\gamma' D\gamma''$, where $\gamma' ,\gamma'' \in
\mathcal{D}$, and $\gamma''$ starts either with a peak or with two
consecutive up steps (if it is not the empty path), we arrive at
the following system (where $F$ is defined as in the previous
proposition):
\begin{equation}\label{2system}
 \begin{cases}
  F=1+xF(1+G+qH)\\
  G=x(1+G+qH)\\
  H=x^2 F(1+G+qH)^2 .
 \end{cases}
\end{equation}

Solving for $F$, we find the following expression:
\begin{displaymath}
F(q,x)=\frac{1-2(1-q)x-\sqrt{1-4x+4x^2 -4qx^2}}{2qx}.
\end{displaymath}


Moreover, the explicit expression of $V(q,x)$ (see again the
previous proposition) can be found in \cite{D,FM2}, and is the
following:
\begin{equation}\label{valleys}
V(q,x)=\frac{1-(1-q)x-\sqrt{1-2(1+q)x+(1-q)^2x^2}}{2qx}.
\end{equation}

We can therefore apply the previous proposition, thus obtaining
formula (\ref{2chains}).\cvd

The integer sequence associated with $SC_2 (x)$ starts
$0,0,0,4,30,168,840,3960,18018,80080,\ldots$. We observe that the
terms of the above sequence divided by 2 yield sequence A002740 of
\cite{S}. In terms of Dyck paths, this sequence gives the sum of
the abscissae of the valleys in all Dyck paths of semilength
$n-1$. It would be nice to have a combinatorial explanation of
this fact.

\bigskip

The results of the present section allow us to compute the Hasse
index of order 2 of Dyck lattices. Recall that in \cite{FM2} it is
shown that the Hasse index of order 1 is asymptotically Boolean.

\begin{prop} The Hasse index of order 2 of the class of Dyck lattices is
asymptotically Boolean.
\end{prop}

\emph{Proof.}\quad Since $|\mathcal{D}_n |={2n\choose n}\frac{1}{n+1}$,
from formula (\ref{2coeff}) we get
\begin{displaymath}
i_2 (\mathcal{D}_n)=\frac{sc_2 (\mathcal{D}_n)}{|\mathcal{D}_n|}=
\frac{(n-1)(n-2)(n+1)}{2(2n-1)}\sim \frac{n^2}{4},
\end{displaymath}
which means precisely that the Hasse index of order 2 is
asymptotically Boolean.\cvd

\section{Saturated chains of length 3}

Setting $h=3$ in (\ref{general}) we obtain a formula for the
enumeration of saturated chains of length 3 of Dyck lattices.
Similarly to what we did in the previous section, we observe that
there are three partitions of the integer 3, namely (1,1,1), (2,1)
and (3). Moreover, the unique ``admissible" triple of skew Ferrers
shapes of area 1 is $({\tiny \yng(1)},{\tiny \yng(1)},{\tiny
\yng(1)})$, whereas there are two pairs of skew Ferrers shapes
whose first component have area 1 and whose second component has
area 2, namely $({\tiny \yng(1)},{\tiny \yng(1,1)})$ and $({\tiny
\yng(1)},{\tiny \yng(2)})$, and there are four skew Ferrers shapes
having area 3, i.e. ${\tiny \yng(1,1,1)}$ , ${\tiny \yng(3)}$ ,
${\tiny \young(:\hfil,\hfil:\hfil)}$ and ${\tiny \yng(2,1)}$ .
Unlike the previous case, now we have two shapes (of area 3) each
of which can be endowed with two different Young tableaux
structures. More precisely, we have to consider the two skew Young
tableaux ${\tiny \young(:3,2:1)}$ , ${\tiny \young(:2,3:1)}$ and
the two (skew) Young tableaux ${\tiny \young(3:1,2)}$ , ${\tiny
\young(3:2,1)}$ . Thus, a direct application of formula
(\ref{general}) leads to the following statement.

\begin{prop} The generating series for the number of saturated
chains of length 3 of Dyck lattices is given by
\begin{align}\label{3sat}
SC_3 (x)=\sum_{n\geq 0}\sum_{\gamma \in \mathcal{D}_n}&(6\cdot \#
(du,du,du)_{\gamma}+3\cdot \# (du,ddu)_{\gamma}\nonumber
\\ &+3\cdot \# (du,duu)_{\gamma}+\# (dddu)_{\gamma}+\# (duuu)_{\gamma}\nonumber
\\ &+2\cdot \# (dduu)_{\gamma}+2\cdot \# (dudu)_{\gamma})\; x^n.
\end{align}
\end{prop}

Our next step will be the evaluation of the unknown quantities
appearing in (\ref{3sat}).

\bigskip

Analogously to the case of saturated chains of length 2, we start
by finding an expression of (\ref{3sat}) better suited for
computation.

\begin{prop} Denote with $A(q,x)$, $B(q,x)$ and $C(q,x)$ the generating series of Dyck
paths where $x$ keeps track of the semilength and $q$ keeps track
of the factors $dduu,dudu$ and $duuu$, respectively. Moreover, let
$V(q,x)$ be defined as in the previous section. Finally, let
$F(y,q,x)$ be the generating series of Dyck paths obtained from
the series $F(q,x)$ defined in the previous section by adding the
indeterminate $y$ keeping track of valleys (i.e. of the factor
$du$). Then
\begin{eqnarray}\label{3chains}
SC_3 (x)&=&2\cdot \left[ \frac{\partial A}{\partial
q}\right]_{q=1}+2\cdot \left[ \frac{\partial B}{\partial
q}\right]_{q=1}+2\cdot \left[ \frac{\partial C}{\partial
q}\right]_{q=1}\nonumber \\
& &+\left[ \frac{\partial^3 V}{\partial q^3}\right]_{q=1}+6\cdot
\left[\frac{\partial^2 F}{\partial y\partial q}-\frac{\partial
F}{\partial q}\right]_{y=q=1}.
\end{eqnarray}
\end{prop}

\emph{Proof.}\quad We start by observing that the knowledge of the
generating series $F(y,q,x)$ allows us to compute the term of
(\ref{3sat}) associated with the pair $(du,duu)$. Indeed, it is
clear that, if we differentiate F with respect to $y$ and $q$ and
then evaluate at $y=q=1$, we obtain the generating series of Dyck
paths with respect to semilength and number of pairs $(du,duu)$.
However, in this way we are going to consider also those pairs in
which the valley $du$ is part of the factor $duu$. Thus, to obtain
what we need, we have to subtract the derivative of $F$ with
respect to $q$, then evaluate at $y=q=1$, which yields the
expression
$$
\left[\frac{\partial^2 F}{\partial y\partial q}-\frac{\partial
F}{\partial q}\right]_{y=q=1}.
$$

Moreover, it is clear that the generating series describing the
distribution of the pair $(du,ddu)$ is the same, and this explains
the coefficient 6 in front of the above displayed expression in
formula (\ref{3chains}).

Finally, the meaning of the partial derivatives of the generating
series $A$, $B$ and $C$ are obvious (notice, in particular, that the
the factors $dddu$ and $duuu$ are clearly equidistributed, so they
are both described by series $C$), as well as the triple partial
derivative of $V$ evaluated in $q=1$, which gives 6 times the
distribution of triples of valleys in Dyck paths.\cvd

\begin{teor} The generating series for the number of saturated
chains of length 3 of $\mathcal{D}_n$ is given by
\begin{equation}\label{3chain}
SC_3 (x)=\sum_{n\geq 0}sc_3 (\mathcal{D}_n )x^n =\frac{P(x) -Q(x)\sqrt{1-4x}}{x(1-4x)^3} ,
\end{equation}
where $P(x)=1-13x+59x^2 -100x^3 +16x^4 +64x^5 =(1-4x)^3 (1-x-x^2
)$ and $Q(x)=1-11x+39x^2 -40x^3 -22x^4$.

The coefficients $sc_3 (\mathcal{D}_n )$ can be expressed as
$$
sc_3 (\mathcal{D}_n )={2n\choose n}\frac{(n^3 -7n+2)(n-2)}{4(n+1)(2n-1)}\qquad (n\geq 2).
$$
\end{teor}

\emph{Proof.}\quad We start by considering the generating series
$F$, $G$, $H$ defined in the previous section. Similarly to what
we did in the above proposition, we need to add an indeterminate
$y$ which will keep track of valleys. Thus, in the following, we
will have $F=F(y,q,x)$, and the same for $G$ and $H$.

Using the same decomposition of Dyck paths described in Theorem
\ref{2teor}, we can now rewrite system (\ref{2system}) taking into
account the presence of the indeterminate $y$, thus obtaining
\begin{equation}\label{3system}
 \begin{cases}
  F=1+xF(1+yG+yqH)\\
  G=x(1+yG+yqH)\\
  H=x^2 F(1+yG+yqH)^2 .
 \end{cases}
\end{equation}
The solution of such a system is the following:
\begin{displaymath}
 \begin{cases}
  F=\frac{1-(1+y-2yq)x-\sqrt{(1+2y+y^2 -4yq)x^2 -2(1+y)x+1}}{2yqx}\\
  G=\frac{(1-(1+y)x-\sqrt{\Delta})(1-(1+y-2yq)x+\sqrt{\Delta})}{4yqx-4y^2 q(1-q)x^2}\\
  H=\frac{(-1+(1+y)x+\sqrt{\Delta})(1-(1+y-2yq)x+\sqrt{\Delta})}{2yqx(4yqx-4y^2 q(1-q)x^2 )} ,
 \end{cases}
\end{displaymath}
where $\Delta=1-2(1+y)x+(1+2y+y^2 -4yq)x^2$.

The expression of $F$ allows us to compute the term of
(\ref{3sat}) associated with the pair $(du,duu)$:
\begin{displaymath}
\left[\frac{\partial^2 F}{\partial y\partial q}-\frac{\partial
F}{\partial q}\right]_{y=q=1}=\frac{-2+15x-30x^2 +10x^3
+(2-11x+12x^2 )\sqrt{1-4x}}{2x(1-4x)\sqrt{1-4x}}.
\end{displaymath}

Recalling the expression of the generating series $V$ reported
in (\ref{valleys}), we obtain:
\begin{displaymath}
\left[ \frac{\partial^3 V}{\partial
q^3}\right]_{q=1}=\frac{3(1-11x+40x^2 -50x^3 +10x^4 -(1-9x+24x^2
-16x^3 )\sqrt{1-4x})}{x(1-4x)^2 \sqrt{1-4x}}.
\end{displaymath}

The generating series $A$ and $B$ can be easily computed starting
from the functional equations they satisfy, which can be found in
\cite{STT} and are reported below for the reader's convenience:
\begin{displaymath}
\begin{cases}
x(q+(1-q)x)A^2 -(1+(1-q)(x-2)x)A+1-(1-q)x=0\\
xB^2 +((1-q)(x-1)x-1)B+(1-q)x+1=0.
\end{cases}
\end{displaymath}

More precisely, we obtain the following expressions:
\begin{equation}
\begin{cases}
A(q,x)=\frac{-1+2(1-q)x-(1-q)x^2 +\sqrt{1-4x+2(1-q)x^2 +(1-q)^2
x^4}}{-2x(q+(1-q)x)}\\
B(q,x)=\frac{1+(1-q)x-(1-q)x^2-\sqrt{1-2(1+q)x-(5-4q-q^2 )x^2
-2(1-q)^2 x^3 +(1-q)^2 x^4}}{2x} .
\end{cases}
\end{equation}

Differentiating with respect to $q$ and evaluating at $q=1$ we
then obtain:
\begin{displaymath}
\begin{cases}
\left[ \frac{\partial A}{\partial q}\right]_{q=1}=\frac{1-5x+5x^2
-(1-3x+x^2 )\sqrt{1-4x}}{2x\sqrt{1-4x}}\\
\left[ \frac{\partial B}{\partial
q}\right]_{q=1}=\frac{1-3x-(1-x)\sqrt{1-4x}}{2\sqrt{1-4x}} .
\end{cases}
\end{displaymath}

Instead the computations related to the generating series $C$ are
a little bit more complicated. Again in \cite{STT} we find the
following functional equation satisfied by $C$:
\begin{displaymath}
qxC^3 +(3(1-q)x-1)C^2 -(3(1-q)x-1)C+(1-q)x=0.
\end{displaymath}

Differentiating both sides with respect to $q$ and then solving
for $\frac{\partial C}{\partial q}$ yields:
\begin{displaymath}
\frac{\partial C}{\partial q}=-\frac{xC^3 -3xC^2 +3xC-x}{3qxC^2
+2(3(1-q)x-1)C-3(1-q)x+1}.
\end{displaymath}

Now, evaluating at $q=1$ and recalling that
$C(1,x)=\frac{1-\sqrt{1-4x}}{2x}$ is the generating series of
Catalan numbers, we get the following:
\begin{displaymath}
\left[ \frac{\partial C}{\partial q}\right]_{q=1}=\frac{-1+6x-9x^2
+2x^3 +(1-4x+3x^2)\sqrt{1-4x}}{x(1-4x-\sqrt{1-4x})}.
\end{displaymath}

We finally have all the information needed to compute $SC_3 (x)$
using (\ref{3chains}), and we obtain formula (\ref{3chain}).
A careful algebraic manipulation of this series yields the
stated expression for the coefficients $sc_3 (\mathcal{D}_n )$.\cvd

The integer sequence $sc_3 (\mathcal{D}_n )$ starts
$0,0,0,2,38,322,2112,12210,65494,$ $334334,\ldots$.
Neither this sequence nor such a sequence divided by 2 appear in \cite{S}.

\begin{prop} The Hasse index of order 3 of the class of Dyck lattices is
asymptotically Boolean.
\end{prop}

\emph{Proof.}\quad Since we have not fully explained
the computations needed to derive the coefficients $sc_3 (\mathcal{D}_n )$,
we will provide a proof independent from the explicit knowledge of such coefficients.

Since series (\ref{3chain}) can be rewritten
as:
\begin{displaymath}
\frac{1}{x}\left( 1-x-x^2 -\frac{Q(x)}{(1-4x)^{5/2}}\right) ,
\end{displaymath}
when $n$ is sufficiently large we have
\begin{eqnarray*}
sc_3 (\mathcal{D}_n )=[x^n ]SC_3 (x)=-[x^{n+1}]Q(x)(1-4x)^{-5/2}.
\end{eqnarray*}

Using Darboux's theorem, we get
\begin{displaymath}
sc_3 (\mathcal{D}_n )\sim -\frac{Q(\xi )}{\xi^{n+1}}\; \frac{(n+1)^{5/2-1}}{\Gamma (5/2)},
\end{displaymath}
where $\xi =\frac{1}{4}$. Since $Q(\xi )=\frac{3}{128}$ and $\Gamma (\frac{5}{2})=\frac{3\sqrt{\pi }}{4}$, we obtain
\begin{displaymath}
sc_3 (\mathcal{D}_n )\sim \frac{2^{2n-3}n^{3/2}}{\sqrt{\pi }}.
\end{displaymath}

Recalling that $|\mathcal{D}_n |\sim \frac{4^n}{n\sqrt{n\pi }}$, we finally have
\begin{displaymath}
i_3 (\mathcal{D}_n )=\frac{sc_3 (\mathcal{D}_n )}{|\mathcal{D}_n
|}\sim \frac{n^3}{8} .
\end{displaymath}
\cvd

\section{Conclusions and further work}

We have derived a general formula for the enumeration of saturated
chains of any fixed length $h$ in Dyck lattices. However, we have
applied such a formula only when $h$ is small (namely $h=2,3$).
When $h$ becomes bigger, computations become much more
complicated. Is it possible to conceive a different approach more
suitable for effective computation?

We have proved that the Hasse indexes of order 1,2 and 3 of Dyck
lattices are asymptotically Boolean. The obvious conjecture is
that the Hasse index of any order $h$ is asymptotically Boolean
too.

Is it possible to extend our approach to enumerate chains in Dyck
lattices?

The problem of enumerating (saturated) chains can also be posed
for other classes of posets. In this context, it would be
interesting to find analogous results in the case of Motzkin and
Schr\"oder lattices.


\begin{thebibliography}{99}
\bibitem[BLL]{BLL} F. Bergeron, G. Labelle, P. Leroux,\quad \emph{Combinatorial Species and Tree-like
Structures},\quad Encyclopedia of Mathematics and Its
Applications, \textbf{67}, Cambridge University Press, Cambridge,
1998.
\bibitem[CPQS]{CPQS} W. Y. C. Chen, S. X. M. Pang, E. X. Y. Qu, R. P. Stanley,\quad \emph{Pairs of
noncrossing free Dyck paths and noncrossing partitions},\quad
Discrete Math.\quad \textbf{309} (2009) 2834--2838.
\bibitem[D]{D} E. Deutsch,\quad \emph{Dyck path enumeration},\quad
Discrete Math.\quad \textbf{204} (1999) 167--202.
\bibitem[FM1]{FM1} L. Ferrari, E. Munarini,\quad \emph{Lattices of paths: representation theory and
valuations},\quad J. Comb.\quad \textbf{2} (2011) 265--292.
\bibitem[FM2]{FM2} L. Ferrari, E. Munarini,\quad \emph{Enumeration of edges in the Hasse diagram of
some lattices of paths},\quad preprint.
\bibitem[FP]{FP} L. Ferrari, R. Pinzani,\quad \emph{Lattices of lattice
paths},\quad J. Statist. Plann. Inference \textbf{135} (2005)
77--92.
\bibitem[F]{F} P. Flajolet,\quad \emph{Polya festoons},\quad INRIA Research Report, N. 1507, September 1991, 6pp.
\bibitem[STT]{STT} A. Sapounakis, I. Tasoulas, P. Tsikouras,\quad \emph{Counting strings in Dyck
paths},\quad Discrete Math.\quad \textbf{307} (2007) 2909--2924.
\bibitem[S]{S} N. J. A. Sloane,\quad \emph{The On-Line Encyclopedia of
Integer Sequences},\quad at http://oeis.org.
\end{thebibliography}
\end{document}